\documentclass[a4,12pt,multicol,makeidx]{article}
\def\origin{%p'=(p_1,...,p_m), 
  \clearpage
\vskip-\baselineskip\vskip-\topskip%
  \vbox to 0pt{\vskip-1in%
    \hbox to 0pt{\hskip-1in%
      \hbox to 0pt{\vrule width 1cm height .4pt depth 0mm\hss}%
      \vbox to 0pt{\hrule width .4pt height 0pt depth 1cm\vss}%
    \hss}%
  \vss}%%
  \vskip-\baselineskip
  \vbox to 0pt{\vskip-1in\vskip3cm%
    \hbox to 0pt{\hskip-1in\hskip3cm%
      \hbox to 0pt{\hss\vrule width 2cm height .4pt depth 0mm\hss}%
      \vbox to 0pt{\vss\hrule width .4pt height 1cm depth 1cm\vss}%
    \hss}%
  \vss}%
\vskip5mm\hskip10mm (3cm,3cm)
}%
\def\a{\alpha}  \def\ou{\overline{u}} \def\oI{\overline{I}} \def\ox{\overline{x}} \def\oy{\overline{y}}
\def\ov{\overline{v}}
 \def\l{\lambda}     \def\e{\varepsilon} 
\def\n{\nabla}    
\def\leq{\underline{<}} 
\def\la{\langle} \def\ra{\rangle} \def\hx{\hat{x}} 
\newenvironment{theorem}{%
\par \bigskip \it}{%
\bigskip \par}

\newenvironment{definition}{%
\par \bigskip \it}{%
\bigskip \par}

\makeindex
\title{Homogenization of a class of integro-differential equations with L{\'e}vy operators. 
}
\author{Mariko Arisawa\\ 
Wolfgang Pauli Institute, University of Wien
\\ Nordbergstrasse 15, A1090 Wien, Austria\\
E-mail: mariko.arisawa@univie.ac.at\\
}
\date{}
\begin{document}
\maketitle
\bigskip
%\pagestyle{plain}
%%%%%Section 1.%%%%%
\section{Introduction.} 

$\quad$ We study the periodic homogenization of 
\begin{equation}\label{first}
	u_{\e}
	 -c(\frac{x}{\e})\int_{z\in {\bf R^N}}  
	[u_{\e}(x+z)
	-u_{\e}(x)- {\bf 1}_{|z|\leq 1}\la \n u_{\e}(x),z \ra]q(z) dz -g(\frac{x}{\e})=0
	\quad  {x\in \Omega},
\end{equation}
\begin{equation}\label{dirichlet}
	u_{\e}(x)=\phi(x)
	\quad x\in {\Omega}^c, 
\end{equation}
where the integral term, the L{\'e}vy operator, has the symmetric density
\begin{equation}\label{density}
	q(z)=\frac{1}{|z|^{N+\a}} \quad z\in {\bf R^N}, \quad \a\in (0,2)\quad \hbox{a constant},
\end{equation}
 $\Omega$ a bounded domain in  ${\bf R^N}$, $c(\cdot)$ and $g(\cdot)$ real valued, periodic, continuous functions in $ {\bf T^N}$, 
 $c(x)>\exists c_0>0$, and $\phi$ a 
continuous function defined in $\Omega^c$. 
 We consider (\ref{first})-(\ref{dirichlet}) in the framework of viscosity solutions for the integro-differential equation (PIDE in short), introduced and studied in 
 A. Sayah \cite{saya}, O. Alvarez and A. Tourin \cite{alvtou}, G. Barles, R. Buckdahn and E. Pardoux \cite{bbp}, H. Pham \cite{pham}, M. Arisawa \cite{ar3}, \cite{ar4}, \cite{ar5}, \cite{ar6}, E. Jacobsen and K. Karlsen \cite{jk} and G. Barles and C. Imbert \cite{bi}.  See  M. Crandall, H. Ishii and P.-L. Lions \cite{users}, too. The comparaison and the existence of solutions have been proved  in the  above works. Recently in \cite{ar6},  the equivalence of several existing  definitions was proved (see Definitions 1.1 and 1.2 in below). The homogenization means to get the unique limit $\lim_{\e\to 0}u_{\e}=\ou$, and to characterize $\ou$ by its effective equation.
 In mathematical finances, (\ref{first}) 
and its evolutionary form are used in the stochastic volatility model with  jump processes (see for example 
R. Cont and P. Tankov \cite{cont}, J.P. Fouque, G. Papanicolaou, and K. Sircar \cite{papa}.) 
We use the formal asymptotic expansion method  introduced by A. Bensoussan, J. L. Lions and G. Papanicolaou \cite{blp} for linear PDEs, and then extended to nonlinear problems 
by P.-L. Lions, G. Papanicolaou and S. Varadhan \cite{lpv} (for first-order PDEs ), L.C. Evans \cite{ev2} (for second-order PDEs), in the framework of viscosity solutions. We shall derive the effective PIDE for $\ou$, rigorously. First, we remind two equivalent definitions of viscosity solutions for a class of PIDEs including (\ref{first}): 
\begin{equation}\label{A}
	A(x,u(x),\n u(x),\n^2 u(x),I[{u}](x))=0
	\quad x\in {\Omega}, 
\end{equation}
where $A(x,u,p,Q,I)$$\in C(\Omega \times {\bf R}\times{\bf R^N} \times {\bf S^N} \times {\bf R} )$, 
 $I[u](x)= \int_{z\in {\bf R^N}}  [u(x+z)-u(x)- {\bf 1}_{|z|\leq 1}\la \n u(x),z \ra]q(z) dz$. 
For an upper (resp. lower) semicontinuous function 
$u\in USC({\bf R^N})$ (resp.  $LSC({\bf R^N})$), $(p,X)\in {\bf R^N}\times {\bf S^N}$ is a sub(resp. super)-differential 
 of $u$ at $x$ : if for any $\delta>0$ there exists $\e>0$ such that 
\begin{equation}\label{pX}
	u(x+z)-u(x)\leq  (\hbox{resp}.\geq)   \la p,z \ra + \frac{1}{2} \la Xz,z \ra + (\hbox{resp}. -)\delta |z|^2\quad \forall |z|\leq \e, 
\end{equation}
 Denote the set of all subdifferentials  (resp. superdifferentials) of $u$ at $x$  $J_{{\bf R^N}}^{2,+}u(x)$ (resp. $J_{{\bf R^N}}^{2,-}u(x)$).  
Set 
$
	I_{\nu,\delta}^{1,+}[u,p,X](x) = \int_{|z|\leq \nu}\frac{1}{2}\la(X+2\delta I)z,z \ra q(z)dz
$ (resp.
$
	I_{\nu,\delta}^{1,-}[u,p,X](x) = \int_{|z|\leq \nu}\frac{1}{2}\la(X-2\delta I)z,z \ra q(z)dz
$), and 
$$
	I_{\nu,\delta}^{2}[u,p,X](x) = \int_{|z|> \nu} [u(x+z)-u(x)- {\bf 1}_{|z|\leq 1}\la p,z \ra] q(z)dz.
$$
We use the following two equivalent definitions (see \cite{ar6}). 
\begin{definition}{\bf Definition 1.1.}  A function $u\in USC({\bf R^N})$ (resp. $LSC({\bf R^N})$) is a viscosity subsolution (resp. supersolution) of (\ref{A}), if for any $\hx\in \Omega$, any $(p,X)\in J_{{\bf R^N}}^{2,+}u(\hx)$ (resp. $J_{{\bf R^N}}^{2,-}u(\hx)$), and 
any pair of numbers $(\e,\delta)$ satisfying  (\ref{pX}), 
the following holds 
$$
	A(\hat{x},u(\hat{x}),p,X, I_{\nu,\delta}^{1,+} (\hbox{resp.} I_{\nu,\delta}^{1,-}) [u,p,X](\hx)
	+ I_{\nu,\delta}^{2}[u,p,X](\hx)) \leq (\hbox{resp.} \geq)  0.
$$
 If $u$ is a subsolution and a supersolution , it is called a viscosity solution.
\end{definition}
\begin{definition}{\bf Definition 1.2.}  A function $u\in USC({\bf R^N})$ (resp. $LSC({\bf R^N})$) is a viscosity subsolution (resp. supersolution) of (\ref{A}),  if for any $\hx\in \Omega$, any $\phi\in C^2({\bf R^N})$ such that $u(\hx)=\phi(\hx)$ and $ u-\phi $ takes a global maximum (resp. minimum) at $\hx$, 
$$
	A(\hat{x},u(\hat{x}),\n \phi(\hat{x}),\n^2 \phi(\hat{x}), I[\phi](\hx)) \leq (\hbox{resp.} \geq) 0.
$$
If $u$ is a subsolution and a supersolution , it is called a viscosity solution.
\end{definition}
We sometimes abbreviate "viscosity" to note a (sub or super) solution. The problem (\ref{first}) was chosen for simplicity to illustrate the method. The various generalizations are possible, namely to the nonlinear problem: 
\begin{equation}\label{extention}
	u_{\e}+H(\frac{x}{\e},\n u_{\e}, \n^2 u_{\e}, I[u_{\e}](x))=0. 
\end{equation}
\section{Formal asymptotic expansions.} 
Let $u_{\e}$ be the solution of (\ref{first}), and assume that 
$$
	u_{\e}(x)=\ou(x)+\e^{\a}v(\frac{x}{\e})+o(\e^{\a}) \quad \forall x\in {\bf R^N}.
$$
 Formally, 
$\n u_{\e}(x)=\n \ou(x) + \e^{\a-1}\n_y v(\frac{x}{\e})$, 
	$\n^2 u_{\e}(x)=\n^2 \ou(x) + \e^{\a-2}\n^2_y v(\frac{x}{\e})$, 
and by introducing them into (\ref{first}), we get 
$$
	\ou
	 -c(\frac{x}{\e})\int_{{\bf R^N}}  
	[\ou(x+z)
	-\ou(x)- {\bf 1}_{|z|\leq 1}\la \n \ou(x),z \ra]q(z) dz \qquad\qquad\qquad\qquad
$$
$$
	 -c(\frac{x}{\e})\int_{{\bf R^N}}  
	\e^{\a}[v(\frac{x+z}{\e})
	-v(\frac{x}{\e})- {\bf 1}_{|z|\leq 1}\la \n_{y} v(\frac{x}{\e}),\frac{z}{\e} \ra]q(z) dz 
	=g(\frac{x}{\e})+o(1).
$$
Put $y=\frac{x}{\e}$, and change the variable to $z'=\frac{z}{\e}$. From (\ref{density}), we have
\begin{equation}\label{cell1}
	\ou
	 -c(y)I[\ou](x)
	 -c(y)I[v](y)
	-g(y)=0.
\end{equation} 
Then,  for each fixed $(x,I)\in \Omega\times {\bf R}$ ($I=I[\ou](x)$ in (\ref{cell1})), find a unique number $d(x,I)$ such that 
 there exists  a periodic solution $v(y)$ of 
\begin{equation}\label{ergodic1}
	d(x,I)
	 -c(y)\int_{{\bf R^N}}  
	[v(y+z)
	-v(y)- {\bf 1}_{|z|\leq 1}\la \n_{y} v(y), z \ra]q dz-g-cI=0, 
\end{equation}
in ${\bf T^N}$. In fact, the existence of $d(x,I)$ (in a  weaker sense) was shown in \cite{ar5} (see Theorem 3.1 in below). 
The effective nonlocal operator is defined as $\oI(x,I)=-d(x,I)$ ($(x,I)\in \Omega \times {\bf R}$), and 
 from (\ref{cell1}), (\ref{ergodic1}), we get: 
\begin{equation}\label{oI1}
	\ou+ \oI(x,I[\ou](x))=0\quad x\in \Omega, 
\end{equation}
the effective equation for $\ou$. 
Later, we justify (\ref{oI1}) by a rigorous argument. \\
\section{The derivation of the  effective equation.}
$\quad$To see the existence of $d(x,I)$ in (\ref{ergodic1}), consider the following 
\begin{equation}\label{ul}\l u_{\l}+ H(y,\n u_{\l})
	 -\int_{{\bf R^N}}  
	[u_{\l}(y+z)-u_{\l}(y)
- {\bf 1}_{|z|\leq 1}\la \n u_{\l}(y),z \ra]q dz 
	-g=0, 
\end{equation} for $y\in {\bf T^N}$,  $\l \in (0,1)$,
 $H$, $g$  real valued functions defined in ${\bf T^N}\times {\bf R^N}$, ${\bf T^N}$, periodic and 
Lipschitz continuous in $y$. \\
\begin{theorem} {\bf Theorem 3.1.(\cite{ar5})} Let $H(y,p)=a(y)|p|$ or $0$, where $a(\cdot)\geq \exists a_0>0$ periodic in ${\bf T^N}$, 
 and consider  (\ref{ul}). The following unique number  $d_g$ exists: 
\begin{equation}\label{dg}
	\lim_{\l \downarrow 0} \l u_{\l}(y)=d_g \quad y{\in}{\bf T^N}, 
\end{equation}
and for any $\rho>0$, there are periodic sub and super solutions  $\underline{u}$ and 
 $\ou$ of 
$$
	d_g+ H(\n \underline{u}(y)) 
	 +I[\underline{u}](y) 
	-g\leq \rho, \quad 
	d_g+ H(\n \overline{u}(y))
	 +I	[\overline{u}](y) 
	-g\geq -\rho \quad y\in {\bf T^N}. 
$$
In particular, if $N=1$ the convergence  (\ref{dg}) is uniform, and for $\rho=0$ there exists $u=\underline{u}=\ou$ which satisfies the above at the 
same time. 
\end{theorem}
We refer the readers to \cite{ar5} (Theorem 6.1) for the proof of the above result. \\
{\bf Remark 3.1.} The convergence  (\ref{dg})  is the ergodic property (see M. Arisawa and P.-L. Lions \cite{al} for the case of PDE). For the case of PIDE, (\ref{dg}) holds in more generality, e.g. for $H$$=H(x,\n u,\n^2 u)$  second-order uniformly elliptic fully nonlinear operator(see \cite{ar5} ). In such a case, the nonlocal homogenization (\ref{extention}) can be solved by the same method in this paper. \\

From Theorem 3.1, for any $(x,I)\in \Omega\times {\bf R}$, there is $\exists^{!}d(x,I)\in {\bf R}$ such that 
 for any $\rho>0$ there exist $\underline{v}$, $\ov$, periodic sub and super solutions of 
$$
	d(x,I)	 +c(y)I[\underline{v}](y)
	-g(y) -c(y)I \leq \rho \quad y\in {\bf T^N},
$$
$$
	d(x,I)	 + c(y)I[\overline{v}](y)
	-g(y) -c(y)I \geq- \rho \quad y\in {\bf T^N}. 
$$
Define  $\oI(x,I)=-d(x,I)$ ($(x,I)\in \Omega \times {\bf R}$).  We remark  the following qualitative property, the degenerate version of which was first stated in \cite{bi}. \\
(Uniform subellipticity) There exists $\theta>0$ such that 
\begin{equation}\label{subelliptic}
	\oI(x,I+I') \leq \oI(x,I) - \theta I'\quad \forall I' >0,
	\quad
	\forall (x,I)\in {\Omega}\times {\bf R}. 
\end{equation}

\begin{theorem}  {\bf Theorem 3.2.$\quad$} The effective integro-differential operator $\oI(x,I)$ is 
 continuous in $\Omega \times {\bf R}$, and is uniformly subelliptic (\ref{subelliptic}) 
with $\theta=c_0$.\\
\end{theorem}

$Proof.$ The proofs  are similar to the PDE's  case in 
 \cite{ev2}. We do not rewrite the proof of the continuity, and mimic that of (\ref{subelliptic}) for the reader's convenience. For $I'>0$,  $I\in {\bf R}$, $\rho>0$,  from Theorem 3.1 we can take $v^I$, $v^{I+I'}$  respectively a sub 
and a super solution of 
\begin{equation}\label{vsub}
	d(x,I)-c(y)I[v^I](y)-g(y)-c(y)I \leq \rho  \quad y\in {\bf T^N}.
\end{equation}
\begin{equation}\label{vsuper}
	d(x,I+I')-c(y)I[v^{I+I'}](y)
	-g(y)-c(y)(I+I') \geq -\rho  \quad y\in {\bf T^N}.
\end{equation}
By adding a constant if necessary, we may asume that $v^{I+I'} < v^{I}$. 
Our goal is to prove 
$\oI(x,I+I') \leq \oI (x,I) - c_0 I'$, $\forall (x,I)\in \Omega\times {\bf R}$. 
Assume the contrary, i.e. there exists a constant $l>0$ such that  
$\oI(x,I+I') \geq \oI (x,I) - c_0 I'  + l$, 
and we shall look for a contradiction. We claim that $v^{I+I'}$ is a viscosity supersolution of 
\begin{equation}\label{vclaim}
	-\oI(x,I)
	 -c(y)I[v^{I+I'}](y)
	-g(y)-c(y)I \geq l-\rho  \quad y\in {\bf T^N}.
\end{equation}
To see this, assume that there exists $\phi\in C^2(\bf R^N)$ such that $v^{I+I'}-\phi$ takes a global 
maximum at a point $y_0\in \Omega$, $v^{I+I'}(y_0)=\phi(y_0)$, and 
$$
	\phi(y_0+z)-\phi(y_0)\geq \la \n \phi (y_0), z\ra + \frac{1}{2} \la (\n^2 \phi (y_0)-2 \delta I)z, z\ra  \quad \forall |z|\leq \nu. 
$$
Since  $v^{I+I'}$ is the supersolution of (\ref{vsuper}), by Definition 1.1, 
$$
	-\oI(x,I+I')
	 -c(y_0)\int_{|z|\leq \nu} \frac{1}{2}   \la (\n^2 \phi (y_0)-2 \delta I)z, z\ra q(z) dz
	-c(y_0)\int_{|z|> \nu}
	[v^{I+I'}(y_0+z)
$$
$$
	-v^{I+I'}(y_0)- {\bf 1}_{|z|\leq 1}\la \n \phi (y_0), z \ra]q(z) dz
	-g(y_0)-c(y_0)(I+I') \geq -\rho.
$$
Then, since $c(y_0)> c_0$ 
$$
	-c(y_0)\int_{|z|\leq \nu} \frac{1}{2}   \la (\n^2 \phi (y_0)-2 \delta I)z, z\ra q(z) dz
	\qquad\qquad\qquad\qquad\qquad\qquad\qquad
$$
$$
	-c(y_0)\int_{|z|> \nu}
	[v^{I+I'}(y_0+z)-v^{I+I'}(y_0)- {\bf 1}_{|z|\leq 1}\la \n \phi (y_0), z \ra]q(z) dz
	-g(y_0)-c(y_0)I
$$
$$
	\geq -c(y_0)\int_{|z|\leq \nu} \frac{1}{2}   \la (\n^2 \phi (y_0)-2 \delta I)z, z\ra q(z) dz
	-c(y_0)\int_{|z|> \nu}
	[v^{I+I'}(y_0+z)-v^{I+I'}(y_0)
$$
$$
	- {\bf 1}_{|z|\leq 1}\la \n \phi (y_0), z \ra]q(z) dz-g(y_0)-c(y_0)(I+I') + c_0 I'  \geq \oI(x,I+I')+ c_0 I' -\rho
$$
$$
	\geq \oI(x,I)-c_0 I' + c_0 I' + l-\rho= \oI(x,I) + l-\rho, 
	\qquad\qquad\qquad\qquad\qquad
$$
and (\ref{vclaim}) is confirmed. For $\l>0$ small enough, from (\ref{vclaim}) we have 
$$
	\l v^{I+I'}(y) 
	-c(y)I[v^{I+I'}](y)
	-g(y)-c(y)I  \geq \oI(x,I)+ l -2 \rho \quad \forall y\in {\bf T^N}, 
$$
while for $\l>0$ small enough, 
$\l v^I  -c(y)I [v^I](y)-g(y)-c(y)I\leq  \oI(x,I) + 2\rho$,  in ${\bf T^N}$. 
 From the comparison (\cite{ar3}, \cite{ar4}), by taking $\rho=\frac{l}{8}$ we get 
$
	\sup_{y\in {\bf T^N}} \l(v^I (y)-v^{I+I'}(y) ) \leq 4 \rho - l\leq -\frac{l}{2}, 
$
which contradicts to $v^{I+I'} < v^{I}$. Thus, $\oI$ is  uniformly subelliptic. \\

$\quad$ Now,  we get the effective equation for $\ou=\lim_{\e\to 0}u_{\e}$: 
\begin{equation}\label{comparison}
	u+ \oI(x,I[u](x)) = 0 \quad x\in \Omega, 
\end{equation}
with (\ref{dirichlet}).  The following comparison result holds.\\
\begin{theorem} {\bf Theorem 3.3.$\quad$}
Let $u\in USC({\bf R^N})$ and $v\in LSC({\bf R^N})$ be respectively a sub and a super solution 
 of (\ref{comparison})- (\ref{dirichlet}). Then, 
$u\leq v$ in $\Omega$.
\end{theorem}
$Proof.\quad$ Since $I$ is uniformly subelliptic (Theorem 3.2), the proof is quite similar to those in 
 \cite{ar3}, \cite{ar4} and \cite{bi} (see \cite{ar5}, too). So, we abbreviate it. \\
\section{The justification of the effective equation.}
$\quad$The main result of this paper is the following. \\
\begin{theorem}  {\bf Theorem 4.1.$\quad$} Let $u_{\e}$ be the solution of (\ref{first}). Then, there exists a unique 
$\lim_{\e \to 0} u_{\e}(x)= \exists \ou(x)$,  which is the solution of (\ref{comparison})-(\ref{dirichlet}). 
\end{theorem}
$Proof.\quad$ Put $u^{\ast}(x)=\lim\sup_{\e\to 0,y\to x} u_{\e}(y)$, $u_{\ast}(x)=\lim\inf_{\e\to 0,y\to x} u_{\e}(y)$. 
 As we shall show in below in Lemma 4.2, $u^{\ast}$, $u_{\ast}$ are respectively a sub and a super solution of (\ref{comparison})- (\ref{dirichlet}). Then, 
from the comparison (Theorem 3.3), 
$u^{\ast}\leq u_{\ast}$, and $u^{\ast}\leq u_{\ast}\leq u^{\ast}$   leads  $\exists !\ou = \lim_{\e \to 0} u_{\e}= u^{\ast}=u_{\ast}$ which is the unique 
solution of (\ref{comparison}). To complete the proof, we need the following.\\
\begin{theorem} {\bf Lemma 4.2.$\quad$}
Let $u_{\e}$ be the solution of (\ref{first}). Then, $u^{\ast}$ and $u_{\ast}$ are respectively a sub and a super solution of (\ref{comparison}).
\end{theorem}
$Proof$ $of$ $Lemma \quad 4.2.$ We show that $u^{\ast}$ is a subsolution of (\ref{comparison}). The proof that 
 $u_{\ast}$ is a supersolution  is shown in parallel, and we abbreviate it. 
 Assume that for $\phi\in C^2(\bf R^N)$, $u^{\ast}-\phi$ takes a global maximum at $\hx\in \Omega$ and $u^{\ast}(\hx)=\phi(\hx)$. As usual (\cite{users}), we may 
assume that  $u^{\ast}-\phi$ takes the global "strict" maximum at $\hx$.  
From Definition 1.2 our goal is to show 
 \begin{equation}\label{goal}
	u^{\ast}(\hx)+\oI (\hx,I[\phi](\hx)) \leq 0. 
\end{equation}
 We use the argument by contradiction. Assume the contrary to (\ref{goal}): 
 \begin{equation}\label{3contra}
	\phi(\hx)+\oI (\hx,I[\phi](\hx)) = 3\gamma>0, 
\end{equation}
for $\gamma>0$. 
Since $\oI$ is continuous, there is $U_{r}(\hx)$$=\{x||x-\hx|<r\}$ such that 
$$
	\phi(x)+\oI (x,I[\phi](x)) \geq \gamma>0 \quad \forall x\in U_{r}(\hx). 
$$
Put $I=I[\phi](\hx)$. By Theorem 3.1, a unique number $d(\hx,I)$ exists, and for any $\rho>0$  there exists 
a periodic continuous fuction  $v(y)$ satisfying 
\begin{equation}\label{5.1approxsup}
	d(\hx,I)-cI -cI[v](y)-g(y)\leq \rho, \quad d(\hx,I)-cI -cI[{v}](y)-g(y)\geq -\rho, 
\end{equation}
in ${\bf T^N}$. 
For $\phi_{\e}(x)=\phi(x)+ \e^{\a}v(\frac{x}{\e})$, 
 (\ref{3contra}) implies that $\phi_{\e}$ is a supersolution of 
\begin{equation}\label{original}
	\phi_{\e}-c(\frac{x}{\e}) I[\phi_{\e}](x)
	-g(\frac{x}{\e})\geq \gamma \quad x\in U_r(\hx), 
\end{equation}
for $r>0$ small enough, i.e. for $\psi\in C^2$ such that $\phi_{\e}-\psi$ attains 
a global minimum at $\ox\in U_r(\hx)$, $(\phi_{\e}-\psi)(\ox)=0$, and we can show (Definition 1.2)
\begin{equation}\label{medium}
	\phi_{\e}(\ox)-c(\frac{\ox}{\e}) \int_{{\bf R^N}}  
	[\psi(\ox+z)
	-\psi(\ox)- {\bf 1}_{|z|\leq 1}\la \n \psi(\ox),z \ra]q(z) dz 
	-g(\frac{\ox}{\e})\geq \gamma. 
\end{equation}
For $h(y)=\frac{1}{\e^{\a}}(\psi-\phi)(\e y)$, $(v-h)(y)$ attains 
a global minimum at $\oy=\frac{\ox}{\e}$, as $\phi_{\e}-\psi$ takes 
the global minimum at $\ox$. Since $v$ is a supersolution of (\ref{5.1approxsup}), 
$$
	d(\hx,I)-c(\oy)I 
	 -c(\oy)I[{h}](\oy)
	-g(\oy)\geq -\rho. 
$$
From the assumption (\ref{3contra}), since $I=I[{\phi}](\hx)$
$$
	\phi(\hx)-c(\oy)I [{\phi}](\hx)
	\qquad\qquad\qquad\qquad\qquad\qquad\qquad\qquad\qquad\qquad\qquad\qquad\qquad\qquad\qquad
$$
$$
	 -c(\oy)\int_{{\bf R^N}}  
	[{h}(\oy+z)
	-{h}(\oy)- {\bf 1}_{|z|\leq 1}\la \n {h}(\oy),z \ra]q(z) dz 
	-g(\oy)\geq 3\gamma -\rho. 
$$
By remarking $h(\frac{x}{\e})=\frac{1}{\e^{\a}}(\psi-\phi)(x)$, $\n_y h(y)={\e^{1-\a}} \n_x (\psi-\phi)(\e y)$, 
by changing the variable $y=\frac{x}{\e}$ to $x$, from (\ref{density}), for $\rho=\gamma$, $r$ small enough, 
we get 
$$
	\phi(\hx)
	 -c(\oy)\int_{{\bf R^N}}  
	[{\psi}(\ox+z)
	-{\psi}(\ox)- {\bf 1}_{|z|\leq 1}\la \n {\psi}(\ox),z \ra]q(z) dz 
	-g(\oy)\geq 2\gamma. 
$$
The claim (\ref{medium}) is shown, that is $\phi_{\e}$ is the supersolution of 
 (\ref{original}). From the comparison (\cite{ar3}, \cite{ar4}, \cite{ar6}, \cite{bi}), 
$(u_{\e}-\phi_{\e})(y)$$\leq \max_{U_r^c(\hx)} (u_{\e}-\phi_{\e}) + \gamma$ for $\forall y\in U_r(\hx)$. 
By letting $\e$ to $0$, $y$ to $\hx$, we have 
$(u^{\ast}-\phi)(\hx)$$\leq \max_{U_r(\hx)^c} (u^{\ast}-\phi) + \gamma$. 
Since $\gamma>0$ is arbitrary
$(u^{\ast}-\phi)(\hx)$$\leq \max_{U_r(\hx)^c} (u^{\ast}-\phi)$.
This contradicts to the assumption that $u^{\ast}-\phi$ takes the global strict maximum at $\hx$. 
Therefore, (\ref{3contra}) is false, and (\ref{goal}) is proved, i.e. $u^{\ast}$ is the subsolution of (\ref{comparison}). As mentioned before, the supersolution property of    $u^{\ast}$ is proved similarly. \\

Since we have proved Lemma 4.2, the proof of Theorem 4.1 is completed.\\

%%%%%%%%%%%%%%%%%%%%%%%%%%%%%%%%%%%%%%%%%%%%%%%%%%%%%%%%%%%%%%%%%%%%%%%%%%%

\end{document}